\documentclass[12pt]{amsart}

\usepackage{amsmath,amssymb,url}
\usepackage{hyperref}
\usepackage{amsthm}
\usepackage{graphicx}
\usepackage{tikz}
\usepackage{tikz-cd}
\usepackage{stmaryrd}
\usepackage{mathrsfs,upgreek}
\usepackage{eucal}
\usepackage{enumerate}
\usepackage{changes}
\usepackage{verbatim}

\newtheorem{Thm}{Theorem}[section]

\newtheorem{Lem}[Thm]{Lemma}
\newtheorem{Cor}[Thm]{Corollary}
\theoremstyle{definition}

\newtheorem{Ex}[Thm]{Example}
\theoremstyle{definition}

\numberwithin{equation}{section}

\newcommand{\op}[1]{\textrm{\upshape #1}}

\newcommand{\join}{\vee}

\newcommand{\meet}{\wedge}

\newcommand{\alg}[1]{{\textbf{\upshape #1}}}  %
\newcommand{\vv}[1]{\mathcal {#1}}

\renewcommand{\a}{\alpha}
\renewcommand{\b}{\beta}
\renewcommand{\d}{\delta}

\newcommand{\g}{\gamma}

\renewcommand{\th}{\theta}

\newcommand{\sse}{\subseteq}

\newcommand{\VV}{{\mathbf V}}   

\newcommand{\ib}{\item[$\bullet$]}

\newcommand{\Con}{\operatorname{\mathbf{Con}}}
\newcommand{\con}{\operatorname{Con}}
\newcommand{\adm}{\operatorname{\mathbf{Adm}
}}

\newcommand{\vuc}[2]{#1_1,\dots,#1_{#2}}

\begin{document}
	\markboth{Paolo Aglian\`o, Stefano Bartali, Stefano Fioravanti}
	{On Freese's technique}
	
	%
	
	%
	
	\title{ON FREESE'S TECHNIQUE}

	\author{Paolo Aglian\`o}
	
	\address{Paolo Aglian\`o,
		DIISM,
		Universit\`a di Siena,
		Via Roma 56,
		53100 Siena,
		Italy}
	\email{\tt agliano@live.com}

	\author{Stefano Bartali}
	
	\address{Stefano Bartali,
		DIISM,
		Universit\`a di Siena,
		Via Roma 56,
		53100 Siena,
		Italy}
	\email{\tt bartali.stefano@gmail.com}

	\author{Stefano Fioravanti}
	
	\address{Stefano Fioravanti,
		DIISM,
		Universit\`a di Siena,
		Via Roma 56,
		53100 Siena,
		Italy}
	\email{\tt stefano.fioravanti66@gmail.com}

	\begin{abstract} In this paper we explore some applications of a certain technique (that we call the \emph{Freese's  technique}), which is a tool for identifying  certain lattices as sublattices of the congruence lattice of a given algebra. In particular we will give sufficient conditions for two family of lattices (called the {\em rods} and the {\em snakes}) to be admissible as sublattices of a variety generated by a given algebra, extending an unpublished result of R. Freese and P. Lipparini.
	\end{abstract}

	\maketitle

	

	\section{Introduction}
	
	The genesis of this paper is our interest in a result of R. Freese and P. Lipparini; the result (Theorem \ref{ralph} below), though unpublished, has been in the open for more than ten years, since the authors lectured about it in several conferences. Our interest was sparkled by the technique used to prove the result; this technique, that we call {\em Freese's technique} (see Section \ref{FreeseTech}), is essentially model theoretic in origin and can be used to identify specific sublattices of congruence lattices of algebras in a variety.
	
	Inspired by the methodology used in the proof we replicated it for some interesting lattices. One of our main results regards the presence of particular sublattices in the congruence variety generated by a nonmodular variety. We find an analogue of  \cite{Freese1995} which, given an algebra whose congruence lattice has a pentagon as in Figure \ref{N5}, characterizes the sublattice of $\con(\alg A(\beta))$ generated using the Freese's technique.
    \begin{Thm}\label{Th1}
		Let $\alg A$ be an algebra with a nonmodular congruence lattice. Then there is an $\alg N_5 \leq \con(\alg A)$ (labelled as in Figure \ref{N5}) such that the sublattice $\alg L$ of $\con(\alg A(\beta))$ generated by $\{\alpha_0, \alpha_1, \gamma_0, \gamma_1, \beta_0\}$ is isomorphic to:
		
		\begin{enumerate}
			\item [(1)] $\alg K$ (Figure \ref{KM2}) if and only if  $\alpha_0 \wedge \gamma_1$ and $\gamma_0 \wedge \alpha_1$ are not comparable;

            \item [(2)] $\alg M_1$ (Figure \ref{M1}) if and only if $\alpha_0 \wedge \gamma_1$ and $\gamma_0 \wedge \alpha_1$ are comparable.
		\end{enumerate}
	\end{Thm}

	Furthermore, with the last result in Section \ref{FreeseTech} we provide a characterization of those pentagons as in Figure \ref{N5} which give $\alg K$ or $\alg M_1$ as a sublattice of $\con(\alg A(\beta))$.
	
	\begin{Thm}\label{Th2}
		Let $\alg A$ be an algebra such that $\con(\alg A)$ has a sublattice as the $\alg N_5$ in Figure \ref{N5} and suppose that $[\gamma, \alpha]$ = $\{\gamma, \alpha\}$;  let $\alg L$ be  the sublattice of $\con(\alg A(\beta))$ generated by $\{\alpha_0, \alpha_1, \gamma_0, \gamma_1, \beta_0\}$. Then $\alg L$ is isomorphic to:
		\begin{enumerate}
			\item $\alg K$ if and only if  $(\beta \circ \gamma \circ \beta) \cap \alpha \not\subseteq \gamma$; \item $\alg M_1$ otherwise.
		\end{enumerate}
	\end{Thm}
	
	In the paper we will produce explicit examples of algebras satisfying either case of Theorem \ref{Th2}.

	In Section \ref{Double} we apply the Freese's technique to lattices called \emph{rods} and \emph{snakes} (Figure \ref{rods}). These families of lattices have been introduced by A. Day in \cite{Day1975} and are of interest since the author proved that these lattices are \emph{splitting} in the variety of modular lattices.
	
	Our aim was to prove a {\em local} version of Corollary \ref{ralphcor}, a straightforward consequence of R. Freese and P. Lipparini result (Theorem \ref{ralph}). This Corollary provides a sufficient condition to have rods and snakes as sublattices of lattices in a variety $\vv V$. We extend this result for a wider class of varieties and in order to do so we first use the Freese's technique to prove the following theorem.
	
	\begin{Thm}\label{M33} Let $\alg A$ be an algebra such that $\alg M_3$, with atoms $\alpha, \b, \gamma$, is a sublattice of  $\con(\alg A)$ consisting of pairwise permuting congruences. Then there is an algebra $\alg A(\g) \in \mathbf{SP(A)}$ such that the sublattice of $\op{Con}(\alg A(\g))$ with universe:
		$$
		\{0_{\op{Con}(\alg A(\g))},\eta_0,\a_0 \meet\a_1,\a_0 \meet \b_1, \a_0,\b_0,\g_0,\d_0\}
		$$
		is isomorphic to $\alg M_{3,3}$ (as in Figure \ref{fig14}) and has pairwise permuting congruences.
	\end{Thm}
	
	This Theorem gives a sufficient condition to generate a lattice $\alg M_{3,3}$ (as in Figure \ref{fig14}) starting from an $\alg M_3 \le \con(\alg A)$. Note that the hypothesis of pairwise permuting congruences in $\alg M_3$ is not so restrictive since, by \cite{KearnesKiss2013}, if $\vv V$ has a \emph{weak difference term}  and $\alg M_3 \leq \con(\alg A)$ for some $\alg A \in \vv V$, then the interval in $\con(\alg A)$  between bottom and top of the $\alg M_3$ consists of pairwise permuting congruences.
	
	Using the Freese's Technique we can obtain the following result about rod and snake lattices.
	
	\begin{Thm}\label{snakesThm} Let $\vv V$ be a variety such that there exists $\alg A \in \vv V$ with $\alg M_3 \le \con(\alg A)$ and all the congruences in the $\alg M_3$ pairwise permute. Then $\alg R_n, \alg S_n \in \adm(\vv V)$ for all $n\ge 1$.
	\end{Thm}
	
	This Theorem gives a partial equivalence between the admission of  $\alg M_3, \alg R_n$, and $\alg S_n$ as a sublattices of congruence lattices of algebras in a variety $\vv V$; the admission of lattices as congruence sublattices is a classical problem, starting from the pioneeristic work of Dedekind \cite{Dedekind1900}.
	A modern account can be found in \cite{KearnesKiss2013}.
	
	\section{Preliminaries and notation}\label{notation}
	
	In this section we recall some of the basic definitions in lattice theory. For other elementary concepts in general algebra (such as lattices, algebras, varieties, etc.) our textbook reference is \cite{BurrisSanka}; for more advanced topics (such as abelian congruences and the commutator of congruences) we refer the reader to \cite{KearnesKiss2013}. For the general theory of Mal'cev conditions and Mal'cev classes there is the classical treatment in \cite{Taylor1973} or the more modern approach in \cite{KearnesKiss2013}.
	
	A lattice $\alg L$ is \emph{modular} if, for $a,b,c \in L$, $c \le a$ implies $a \meet (b \join c) \le (a \meet b) \join c$; a lattice is \emph{distributive} if
	$\join$ and $\meet $ distribute w.r.t each other; a lattice is \emph{meet semidistributive} if for $a,b,c \in L$
	$$
	a \meet b = a\meet c \qquad\text{implies}\qquad a \meet b = a \meet (b \join c),
	$$
	is \emph{join semidistributive} if it satisfies the dual property and it is \emph{semidistributive} if it is both join and meet semidistributive. It is clear that a distributive lattice is  both modular and semidistributive; however $\alg N_5$ (Figure \ref{N5M3}) is not modular but it is semidistributive, while  $\alg M_3$ (Figure \ref{N5M3}) is modular but neither join nor meet semidistributive.
	
	\begin{figure}[hbtp]
		\begin{center}
			\begin{tikzpicture}
				\draw (0,0)-- (-1,1) -- (0,2) -- (0.6,1.4) -- (0.6,0.6) -- (0,0);
				\draw[fill] (0,0) circle [radius=0.05];
				\draw[fill] (-1,1) circle [radius=0.05];
				\draw[fill] (0,2) circle [radius=0.05];
				\draw[fill] (0.6,1.4) circle [radius=0.05];
				\draw[fill] (0.6,0.6) circle [radius=0.05];
				
				\draw (5,0) -- (4,1) -- (5,2) -- (5,0) -- (6,1) -- (5,2);
				\draw[fill] (5,0) circle [radius=0.05];
				\draw[fill] (4,1) circle [radius=0.05];
				\draw[fill] (5,2) circle [radius=0.05];
				\draw[fill] (5,1) circle [radius=0.05];
				\draw[fill] (6,1) circle [radius=0.05];
				\node [right] at (1,0) {\footnotesize $\alg N_5$};
				\node [right] at (6,0) {\footnotesize $\alg M_3$};
			\end{tikzpicture}
		\end{center}
		\caption{$\alg N_5$ and $\alg M_3$}\label{N5M3}
	\end{figure}
	
	A lattice $\alg L$ is \emph{projective for a class of lattices} $\vv K$ if for all $\alg M \in \vv K$ and for any onto homomorphism $g: \alg M \rightarrow \alg L$ there is a monomorphism $f: \alg L \rightarrow \alg M$ such that $gf$ is the identity mapping. If $\vv K=\vv L$, the variety of all lattices, we simply say that $\alg L$ is \emph{projective}. The \emph{Whitman's condition} (W) is one of the four properties discovered by P. Whitman \cite{Whitman1941}, characterizing free lattices; a lattice $\alg L$ satisfies (W) if, for all $a,b,c,d \in L$, $a \meet b \le c \join d$ implies either $a \le c \join d$ or $b \le c \join d$ or $a\meet b \le c$ or $a \meet b\le d$.
	
	\begin{Thm} \cite[Theorem $1.4$]{Nation1982} A finite lattice is projective if and only if it is semidistributive and satisfies (W).
	\end{Thm}
	
	A subdirectly irreducible lattice $\alg L$ is \emph{splitting} in a variety $\vv V$ of lattices if  there is a variety $\vv W_\alg L$ such that $\alg L \notin \vv W_\alg L$ and  for any variety $\vv U \sse \vv V$  either $\vv U \sse \vv W_\alg L$ or $\alg L \in \vv U$. It can be shown \cite{McKenzie1972} that in this case $\vv W_\alg L$ is axiomatized by a single equation, called the \emph{splitting equation} of $\alg L$. We can observe that if $\alg L$ is splitting in $\vv V$, then it is splitting in any subvariety of $\vv V$ (hence if $\alg L$ is not splitting in $\vv V$, then it is not splitting in any supervariety of $\vv V$). It is not hard to check that if $\alg L$ is subdirectly irreducible and projective in $\vv V$, then it is splitting in $\vv V$; the converse however fails to hold.

	Finally for any variety $\vv V$,
	$$
	\adm(\vv V) =\{\alg L \mid \alg L \leq \con(\alg A),  \alg A \in \vv  V\},
	$$
	and we say that $\vv V$ \emph{admits} $\alg L$ while otherwise we say that $\vv V$ \emph{omits} $\alg L$.  The \emph{congruence variety} of $\vv V$ is the variety generated by $\adm(\vv V)$ and it is denoted by $\Con(\vv V)$; a variety of lattices is a congruence variety if it is equal to $\Con(\vv V)$ for some variety $\vv V$. In general $\adm(\vv V) \subsetneq \Con(\vv V)$.

	\section{Freese's  technique}\label{FreeseTech}

	Let $\alg A$ be any algebra and $\a \in \con(\mathbf{A})$; by $\alg A^n(\a)$ we denote the subalgebra of $\alg A^n$ whose universe is
	$$
	A^n(\a) :=\{(\vuc an) \mid (a_i,a_j) \in \a, i,j \le n\}.
	$$
	Let $\alg L$ be any lattice, $\vv V$ a variety of algebras, and $\alg A\in \vv V$; we would like to determine under which conditions on the congruence lattice of $\alg A$
	we can infer that $\alg L \in \Con(\vv V)$. A stronger question is:
	when is there  an $\a \in \con(\mathbf{A})$ such that $\alg L \le \op{Con}(\alg A^n(\a))$?
	The idea of looking for sublattices of $\op{Con}(\alg A^n(\a))$ is essentially due to R. Freese, who is the author or the coauthor of almost all the printed material on the subject. For this reason through all the paper we will refer to this as \emph{Freese's  technique}. If $n=2$, we drop the superscript and $\alg A(\a) = \a$. In this case we refer to the Freese's  technique as \emph{duplication technique}, briefly described in \cite[pages $96-101$]{FreeseMcKenzieMcNultyTaylor}.
	
	This is a list of results that can be (and have been) proved using the Freese's  technique; here $\vv M_p$ is the congruence variety of vector spaces over a field of characteristic $p$.
	
	\begin{enumerate}
		\ib There are nonmodular varieties of lattices that are not congruence varieties (Nation \cite{Nation1974} without Freese's  technique, then Freese, see \cite{Freese1995}).
		\ib Every modular congruence variety $\vv V$ consists entirely of {\em arguesian} lattices (Freese and J\'onsson,\cite{FreeseJonsson1976}), hence the variety of all modular lattices is not a congruence variety.
		\ib Every nonmodular congruence variety contains $\Con(\vv P)$, where $\vv P$ is Polin's variety (Freese and Day, \cite{DayFreese1980}).
		\ib Every modular nondistributive congruence variety $\vv V$ contains $\vv M_p$, for some $p$ that is prime or $0$ (Freese, Herrman and Huhn, \cite{FreeseHerrmanHuhn1981}).
		\ib If $\vv V$ is a variety that is not congruence meet semidistributive, then $\Con(\vv V)$ contains $\vv M_p$, for some $p$ that is prime or $0$ (Freese and Lipparini, unpublished).
	\end{enumerate}
	
	We now take a closer look at the duplication technique, that was first introduced in \cite{FreeseJonsson1976}. If $\alg A$ is any algebra and $\a \in \con(\alg A)$, then $\alg A(\a) \le \alg A^2$ and in particular it is a subdirect product of two copies of $\alg A$. We can use the following notation for congruences on subdirect products: if $\th \in \con(\alg A)$, then
	\begin{enumerate}
		\ib $\th_0 = \{( (a_0,a_1),(b_0,b_1)) \mid (a_0,b_0) \in \th\}$;
		\ib $\th_1 = \{( (a_0,a_1),(b_0,b_1)) \mid (a_1,b_1) \in \th\}$;
		\ib $\eta_0$,$\eta_1$ are the kernels of the projections of $\alg A(\a)$ over the two factors.
	\end{enumerate}
	By the Homomorphism  Theorem  \cite[Theorem $6.11$]{BurrisSanka} $\op{Con} (\alg A(\a)/\eta_0) \cong \op{Con} (\alg A(\a)/\eta_1) \cong \con (\alg A)$; moreover the following facts are easy to prove  (see also \cite{FreeseJonsson1976}).
	
	
	\begin{Lem}\label{doublinglemma}  Let $\alg A$ be an algebra, $\g,\th \in \con(\alg A)$. Then the following holds in $\op{Con}(\alg A(\g))$:
		\begin{enumerate}
			\item if $\g \le \th$, then $\th_0=\th_1$;
			\item $\th_i = \eta_i \join (\th_0 \meet \th_1)$;
			\item $\g_0 = \eta_0 \join \eta_1$.
		\end{enumerate}
	\end{Lem}

	\begin{figure}[hbtp]
		\begin{center}
			\begin{tikzpicture}
				\draw (0,0)-- (-1,1) -- (0,2) -- (0.6,1.4) -- (0.6,0.6) -- (0,0);
				\draw[fill] (0,0) circle [radius=0.05];
				\draw[fill] (-1,1) circle [radius=0.05];
				\draw[fill] (0,2) circle [radius=0.05];
				\draw[fill] (0.6,1.4) circle [radius=0.05];
				\draw[fill] (0.6,0.6) circle [radius=0.05];
				\node [right] at (0,-.1) {\footnotesize $0_\alg A$};
				\node [right] at (0.6,0.6) {\footnotesize $\g$};
				\node [right] at (0.6,1.4) {\footnotesize $\a$};
				\node [left] at (-1,1) {\footnotesize $\beta$};
				\node [right] at (0,2.1) {\footnotesize $\delta$};
			\end{tikzpicture}
		\end{center}
		\caption{$\alg N_5$}\label{N5}
	\end{figure}

	Thus suppose that $\alg A$ is an algebra that is not congruence modular; then $\alg N_5 \le \con(\alg A)$. If $\alg N_5$ is as in Figure \ref{N5}, let us apply the duplication technique to $\g$; then Lemma \ref{doublinglemma} tells us that inside $\op{Con}(\alg A(\g))$ there is the sublattice in Figure \ref{L15lat}.
	
	\begin{figure}[h]
		\begin{center}
			\begin{tikzpicture}
				\draw (0,0)-- (-1.5,1.5) -- (0,3) -- (1.5,1.5) -- (0,0) -- (0,1.5) -- (-1.5,3) -- (0,4.5) -- (1.5,3)-- (0,1.5);
				\draw (-1.5,1.5) -- (-1.5,3);
				\draw (0,3) -- (0,4.5);
				\draw (1.5,1.5) --(1.5,3);
				\draw[fill] (0,0) circle [radius=0.05];
				\draw[fill] (-1.5,1.5) circle [radius=0.05];
				\draw[fill] (0,3) circle [radius=0.05];
				\draw[fill] (1.5,1.5) circle [radius=0.05];
				\draw[fill] (0,1.5) circle [radius=0.05];
				\draw[fill] (-1.5,3) circle [radius=0.05];
				\draw[fill] (0,4.5) circle [radius=0.05];
				\draw[fill] (1.5,3) circle [radius=0.05];
				\draw[fill] (0,3.75) circle [radius=0.05];
				\node [right] at (0,-.1) {\footnotesize $\eta_0 \meet \eta_1$};
				\node [left] at (-1.5,1.5) {\footnotesize $\eta_0$};
				\node [right] at (1.5,1.5) {\footnotesize $\eta_1$};
				\node [left] at (-1.5,3) {\footnotesize $\beta_0$};
				\node [right] at (1.5,3) {\footnotesize $\beta_1$};
				\node [right] at (0,3) {\footnotesize $\g_0=\g_1$};
				\node [right] at (0,1.5) {\footnotesize $\beta_0 \meet \beta_1$};
				\node [right] at (0,3.75) {\footnotesize $\a_0=\a_1$};
				\node [right] at (0,4.5) {\footnotesize $\d_0=\d_1$};
			\end{tikzpicture}
		\end{center}
		\caption{$\alg L_{14}$}\label{L15lat}
	\end{figure}

	It follows that if $\vv V$ is not congruence modular then $\alg L_{14} \in \adm(\vv V)$. Now $\VV(\alg N_5)$ is clearly a subvariety of $\VV(\alg L_{14})$; it is also proper, since by J\'onsson Theorem \cite{Jonsson1967}, $\alg L_{14} \notin \VV(\alg N_5)$. So any variety $\vv U$
	of lattices such that $\alg N_5 \in \vv U$ but $\alg L_{14} \notin \vv U$ cannot be a congruence variety. Hence there are nonmodular varieties of lattices that are not congruence varieties. Note that $\alg L_{14}$ is isomorphic to the sublattice generated by the two copies of $\alg N_5$ inside $\alg A(\g)$ and it is thus uniquely determined. R. Freese in \cite{Freese1995} observed that if we try to apply the same duplication construction to $\a$ or $\b$ then the sublattice generated by the two pentagons is no longer uniquely determined.
	
	Let $\alg M_1$ be the lattice in Figure \ref{M1}. In the following Lemma we will discuss the case in which we apply the duplication construction to $\b$.
	
	\begin{Lem}\label{lemma1} Let $\alg A$ be an algebra such that $\con(\mathbf{A})$ has a sublattice as the $\alg N_5$ in Figure \ref{N5}. Then the following holds in $\op{Con}(\alg A(\b))$:
		\begin{enumerate}
			\item $\g_0 \meet \g_1 \le \a_0 \meet \g_1, \g_0 \meet \a_1, \a_0 \meet \a_1$;
			\item $\a_0 \meet \a_1 \not\le \a_0 \meet \g_1, \g_0 \meet \a_1, \g_0 \meet \g_1$;
			\item  if $\a_0 \meet \g_1$ and $\g_0 \meet \a_1$ are comparable congruences, then
			$$
			\a_0 \meet \g_1 =\g_0 \meet \a_1= \g_0\meet\g_1;
			$$
			\item$\g_0 \join (\a_0 \meet \a_1) = \a_0$ and $\g_1 \join (\a_0 \meet \a_1) = \a_1.$
		\end{enumerate}
	\end{Lem}
	All verifications are routine; however it follows that:
	
	\begin{Cor}\label{CorM1}
		Let $\alg A$ be an algebra such that $\con(\alg A) $ has a sublattice as the $\alg N_5$ in Figure \ref{N5}. Let $\a_0 \meet \g_1$ and $\g_0 \meet \a_1$ be comparable congruences in $\op{Con}(\alg A(\b))$.
		Then the sublattice of $\op{Con}(\alg A(\b))$ generated by $\{\alpha_0, \alpha_1, \gamma_0 ,\gamma_1, \beta_0\}$ is isomorphic to $\alg M_1$.
	\end{Cor}
	\begin{proof}If the hypotheses hold, then Lemmata \ref{doublinglemma} and \ref{lemma1} guarantee that all the meets and some of the joins are exactly as in Figure \ref{M1}. Indeed, first we can observe that $\op{Con}(\alg A(\g_1)/\eta_i)$ is isomorphic to $\con(\alg A)$ via the natural isomorphism $f(\lambda) = \lambda_i$, for $i =1,2$. Thus the joins and meets  in Figure \ref{M1} of elements above $\eta_{0}$ and $\eta_{1}$ are correctly depicted. Furthermore, by Lemma \ref{lemma1}, $\a_0 \wedge \gamma_1 = \a_1 \wedge \gamma_0 = \g_0 \wedge \g_1$. We note that $\b_0 \wedge  \a_0 \wedge \a_1 = \b_0 \wedge  \a_0 \wedge \beta_0 \wedge \a_1 = \eta_0 \wedge \eta_1 = 0_{\alg A(\beta)}$ and $\a_0 \wedge \eta_1 = \a_0 \wedge \b_0  \wedge \eta_1 = \eta_0 \wedge \eta_1$. Hence, all the meets in Figure \ref{M1} are correctly displayed. For the joins we can see that $\g_0 \join \eta_1 = \g_0 \join \eta_0 \join \eta_1 = \g_0 \join \b_0 = \d_0$. From this equations, their symmetric versions obtained exchanging $1$ and $0$, and trivial inclusions, we conclude that the joins of the congruences in the set $\{\eta_0, \eta_1, \g_0, \g_1, \a_0, \a_1, \b_0, \d_0\}$ are correctly depicted. Moreover, by Lemma \ref{doublinglemma}, $\b_0 \vee (\g_0 \wedge \g_1) = \b_0 \vee \eta_ 0 \vee (\g_0 \wedge \g_1) = \b_0 \vee \g_0 = \d_0$ and $\g_0 \vee (\a_0 \vee \a_1) = \eta_0 \vee (\a_0 \vee \a_1) = \a_0$. The correctness of the remaining joins follows from immediate inclusions or by symmetry exchanging $1$ and $0$.

	\end{proof}
	
	\begin{figure}[htbp]
		\begin{center}
			\begin{tikzpicture}
				\draw (0,0) -- (-2,1.5) -- (-2,3) -- (0,4.5) -- (0,3) -- (0,4.5) -- (2,3) -- (2,1.5) -- (0,0);
				\draw (0,0) -- (0,1.5);
				\draw (-2,1.5) -- (0,3) -- (2,1.5) -- (0,0);
				\draw (-2,3) -- (0,1.5) -- (2,3);
				\draw (-2,2.2) -- (0,0.7) -- (2,2.2);
				\draw[fill] (0,0) circle [radius=0.05];
				\node [right] at (0,-.2) {\footnotesize $0_{\alg A(\b)}$};
				\draw[fill] (-2,1.5) circle [radius=0.05];
				\node [left] at (-2,1.5) {\footnotesize $\eta_0$};
				\draw[fill] (-2,3) circle [radius=0.05];
				\node [left] at (-2,3) {\footnotesize $\a_0$};
				\draw[fill] (0,4.5) circle [radius=0.05];
				\node [right] at (0,4.5) {\footnotesize $\d_0=\d_1$};
				\draw[fill] (0,3) circle [radius=0.05];
				\node [right] at (0,3) {\footnotesize $\b_0=\b_1$};
				\draw[fill] (2,3) circle [radius=0.05];
				\node [right] at (2,3) {\footnotesize $\a_1$};
				\draw[fill] (2,1.5) circle [radius=0.05];
				\node [right] at (2,1.4) {\footnotesize $\eta_1$};
				\draw[fill] (0,1.5) circle [radius=0.05];
				\node [right] at (0,1.5) {\footnotesize $\a_0\meet \a_1$};
				\draw[fill] (-2,2.2) circle [radius=0.05];
				\node [left] at (-2,2.2) {\footnotesize $\g_0$};
				\draw[fill] (0,0.7) circle [radius=0.05];
				\node [right] at (0,0.7) {\footnotesize $\g_0\meet \g_1$};
				\draw[fill] (2,2.2) circle [radius=0.05];
				\node [right] at (2,2.2) {\footnotesize $\g_1$};
			\end{tikzpicture}
		\end{center}
		\caption{$\alg M_1$}\label{M1}
	\end{figure}
	
	We have seen what happens if $\a_0 \meet \g_1$ and $\g_0 \meet \a_1$ are comparable. The case $\a_0 \meet \g_1$ and $\g_0 \meet \a_1$ not comparable will be characterized  with the next Lemmata.
	
	\begin{Lem}\label{lemma2} Let $\alg A$ be an algebra such that $\con(\alg A)$ has a sublattice as the $\alg N_5$ in Figure \ref{N5}. Assuming that $\a_0 \meet \g_1$ and $\g_0 \meet \a_1$ are incomparable in $\con(\alg A(\b))$, we have:
		$$
		(\a_0 \meet \g_1) \join (\g_0 \meet \a_1) < \a_0 \meet \a_1.
		$$
	\end{Lem}
	\begin{proof} Let $\alpha, \gamma$ be as in the hypothesis. We will show that for any $(a,b) \in \a\setminus \g$,
		$$
		((a,a),(b,b)) \notin (\a_0 \meet \g_1) \join (\g_0 \meet \a_1);
		$$
		this is enough to prove the thesis.
		
		Suppose that $((a,a), (b,b)) \in (\a_0 \meet \g_1) \join (\g_0 \meet \a_1)$; then there is an $n \in \mathbb N$ and $(u_0,v_0), \dots,(u_n,v_n) \in \b$ such that $(a,a) = (u_0,v_0)$, $(b,b)  = (u_n, v_n)$ and
		\begin{align*}
			&(u_i,v_i) \mathrel{\a_0 \meet \g_1} (u_{i+1},v_{i+1}) \qquad\text{for $i$ even}\\
			&(u_i,v_i) \mathrel{\g_0 \meet \a_1} ( u_{i+1},v_{i+1}) \qquad\text{for $i$ odd}.
		\end{align*}
		Thus we have the following relations
		\begin{center}
			\begin{tikzcd}
				a \arrow[dd,dash]\arrow[r,"\a", dash]&u_1\arrow[dd,"\b",dash]\arrow[r,"\g", dash]&u_2\arrow[dd,"\b",dash]\arrow[r,"\a" dash]&u_3\arrow[dd,"\b",dash]\arrow[r,dash,dotted]&\arrow[r,dash,dotted]&u_{n-1}\arrow[dd,"\b",dash]\arrow[r,"\a", dash]&b\arrow[dd,dash] \\
				\\
				a\arrow [r,"\g",dash]&v_1\arrow[r,"\a", dash]&v_2\arrow[r,"\g", dash]&v_3\arrow[r,dash,dotted]&\arrow[r,dash,dotted]&v_{n-1}\arrow[r,"\g", dash]&b
			\end{tikzcd}
		\end{center}
		Since $\g \le \a$ we get $(u_1,v_1) \in \a$ and thus $(u_1,v_1)  \in \a \meet \b = 0_\alg A$, hence $u_1=v_1$. Continuing this argument we get $u_i=v_i$ for $i \le n$, so
		$$
		a \mathrel{\g} v_1=u_1 \mathrel{\g} u_2=v_2 \mathrel{\g} \dots \mathrel{\g} b.
		$$
		Hence $(a,b) \in \g$ that is a contradiction. From this we can conclude that $((a,a),( b,b)) \notin (\a_0 \meet \g_1) \join (\g_0 \meet \a_1)$ and thus $(\a_0 \meet \g_1) \join (\g_0 \meet \a_1) < \a_0 \meet \a_1$.
	\end{proof}
	
	\begin{figure}[htbp]
		\begin{center}
			\begin{tikzpicture}
				\draw (0,0) -- (-2,1.5) -- (-2,3.5) -- (0,5) -- (0,3) -- (0,5.) -- (2,3.5) -- (2,1.5) -- (0,0);
				\draw (0,0) -- (0,0.7);
				\draw (-2,1.5) -- (0,3) -- (2,1.5) -- (0,0);
				\draw (-2,3.5) -- (0,2)-- (0,1.5) -- (0,2) -- (2,3.5);
				\draw (-2,2.2) -- (0,0.7) -- (2,2.2);
				\draw (0,1.5) -- (-.48,1.05);
				\draw (0,1.5) -- (.48,1.05);
				\draw[fill] (0,0) circle [radius=0.05];
				\node [right] at (0,-.2) {\footnotesize $0_{\alg A(\b)}$};
				\draw[fill] (-2,1.5) circle [radius=0.05];
				\node [left] at (-2,1.5) {\footnotesize $\eta_0$};
				\draw[fill] (-2,3.5) circle [radius=0.05];
				\node [left] at (-2,3.5) {\footnotesize $\a_0$};
				\draw[fill] (0,5) circle [radius=0.05];
				\node [right] at (0,5) {\footnotesize $\d_0=\d_1$};
				\draw[fill] (0,3) circle [radius=0.05];
				\node [right] at (0,3) {\footnotesize $\b_0=\b_1$};
				\draw[fill] (2,3.5) circle [radius=0.05];
				\node [right] at (2,3.5) {\footnotesize $\a_1$};
				\draw[fill] (2,1.5) circle [radius=0.05];
				\node [right] at (2,1.4) {\footnotesize $\eta_1$};
				\draw[fill] (0,1.5) circle [radius=0.05];
				\node [right] at (0,1.5) {\footnotesize $\th$ };
				\draw[fill] (-2,2.2) circle [radius=0.05];
				\node [left] at (-2,2.2) {\footnotesize $\g_0$};
				\draw[fill] (0,0.7) circle [radius=0.05];
				\node [right] at (0,0.7) {\footnotesize $\g_0\meet \g_1$};
				\draw[fill] (2,2.2) circle [radius=0.05];
				\node [right] at (2,2.2) {\footnotesize $\g_1$};
				\draw[fill] (0,2) circle [radius=0.05];
				\node [right] at (0,2) {\footnotesize $\a_0\meet \a_1$};
				\draw[fill] (-.48,1.05) circle [radius=0.05];
				\node [left] at (-.48,1.05) {\footnotesize $\g_0 \meet \a_1$};
				\draw[fill] (.48,1.05) circle [radius=0.05];
				\node [right] at (.48,1.05) {\footnotesize $\a_0\meet \g_1$};
			\end{tikzpicture}
		\end{center}
		\caption{$\alg K$}\label{KM2}
	\end{figure}
	We are ready to present a key theorem for this section.

\begin{Thm}\label{Th1-2}
		Let $\alg A$ be an algebra such that $\con(\alg A)$ has a sublattice as the $\alg N_5$ in Figure \ref{N5} and suppose that $[\gamma, \alpha]$ = $\{\gamma, \alpha\}$. Let $\alg L$ be the sublattice of $\con(\alg A(\beta))$ generated by $\{\alpha_0, \alpha_1, \gamma_0, \gamma_1, \beta_0\}$. Then:
		
		\begin{enumerate}
			\item [(1)] $\alpha_0 \wedge \gamma_1$ and $\gamma_0 \wedge \alpha_1$ are not comparable if and only if $\alg L$ is isomorphic to $\alg K$;

            \item [(2)] $\alpha_0 \wedge \gamma_1$ and $\gamma_0 \wedge \alpha_1$ are comparable if and only if $\alg L$ is isomorphic to $\alg M_1$.
		\end{enumerate}
		
	\end{Thm}
	
	\begin{proof}
        For $(\Rightarrow)$ in $(1)$ we define $\theta = (\alpha_0 \wedge \gamma_1) \vee (\gamma_0 \wedge \alpha_1)$ and, by Lemma \ref{lemma2}, $\theta < \a_0 \meet \a_1$. We can see that $(\alpha_0 \wedge \gamma_1) \wedge (\gamma_0 \wedge \alpha_1) = \gamma_0 \wedge \gamma_1$ is a consequence of the definition of the four congruences involved in the equation.
		Furthermore, we can observe that $\alpha_i$ covers $\gamma_i$, since $\con(\alg A) \cong  \con(\alg A(\beta)/\eta_i)$ through the natural isomorphism $f(\lambda) = \lambda_i$. This implies that $\theta \join \g_0 = \alpha_0$ and symmetrically  $\theta \join \g_1 = \alpha_1$. The other meets and joins can be verified using Lemma \ref{lemma1}.  Hence, $\alg K$ is isomorphic to the sublattice of $\con(\alg A(\beta))$ generated by $\{\alpha_0, \alpha_1, \gamma_0,$ $\gamma_1, \beta_0\}$.

		The proof  of $(\Rightarrow)$ in $(2)$ follows from Corollary \ref{CorM1}, which ensures that $\alg M_1$ is isomorphic to the sublattice of $\con(\alg A(\beta))$ generated by $\{\alpha_0, \alpha_1, \gamma_0,$ $\gamma_1, \beta_0\}$.
		
		Thus, the two claims in $(1)$ and $(2)$ hold since ($\Leftarrow$) of both are trivial.
	\end{proof}
    Using the previous theorem we can almost directly prove one of the main results of the paper.
    \begin{proof}[Proof of Theorem \ref{Th1}.]
        By \cite{Dedekind1900}, $\alg N_5 \leq \con(\alg A)$ and, by \cite[2.2]{CrawleyDilworth1973}, every congruence lattice is \emph{weakly atomic} and thus if $\con(\alg A)$ has a sublattice as the $\alg N_5$ in Figure \ref{N5}, 
	   then there exist  $\gamma', \alpha' \in [\gamma, \alpha]$ such that $[\gamma', \alpha']$ = $\{\gamma', \alpha'\}$. Hence, we can observe that the congruences $\{0_{\alg A}, \a', \g', \b, \d\}$ form a pentagon in $\con(\alg A)$ with $[\gamma', \alpha']$ = $\{\gamma', \alpha'\}$ and the claim directly follows from Theorem \ref{Th1-2}.
    \end{proof}
	
	
	The lattice  $\alg K$ in Figure \ref{KM2} is subdirectly irreducible; it is also projective in the variety of all lattices (since it is semidistributive and satisfied the Whitman condition), so it is splitting. The splitting equation of $\alg K$ is a 5-variables equation $p(x,y_0,y_1,z_0,z_1) \le q(x,y_0,y_1,z_0,z_1)$ (see \cite{Freese1995} for the complete expression). R. Freese in \cite{Freese1995} showed that $\alg K$ belongs to any nonmodular congruence variety by proving that the splitting equation of $\alg K$ fails in $\op{Con}(\alg A(\b))$, for some $\alg A$ in the variety and some congruence $\b$ of $\alg A$.
	However now we can extract more information:
	
	\begin{Thm}\label{omitd1} Let $\vv V$ be a variety satisfying a nontrivial idempotent Mal'cev condition and let $\alg A \in \vv V$ be an algebra with a nonmodular congruence lattice. Then there exists an $\alg N_5 \leq \con(\alg A)$ (labelled as in Figure \ref{N5}) such that
    $\alg K \le \op{Con}(\alg A(\b))$.
	\end{Thm}
	\begin{proof} By Theorem \ref{Th1} we have a pentagon in $\con(\alg A)$ such that either $\alg M_1$ or $\alg K$ is a sublattice of $\op{Con}(\alg A(\b))$. Since $\vv V$ satisfies a nontrivial idempotent Mal'cev condition, then $\vv V$ omits $\alg D_1$ as a sublattice of the congruence lattice of its algebras (see \cite[Theorem $4.16$]{KearnesKiss2013}). We can observe that the filter of $\alg M_1$ generated by $\g_0 \wedge \g_1$ is a sublattice isomorphic to $\alg D_1$ (Figure \ref{M1}), thus it cannot appear as sublattice of $\op{Con}(\alg A(\b))$, so $\alg K$ is the only possibility left.
	\end{proof}
	
	Using a strategy similar to Lemma \ref{lemma2} we can also characterize when a configuration as the $\mathbf{N}_5$ in Figure \ref{N5} generates either $\alg M_1$ or $\alg K$ in $\con(\alg A(\beta))$. Let us use the standard notation $\alpha \circ^n \beta = \alpha \circ \cdots \circ \beta$ for the relation product with $n$ factors.
	
	\begin{proof}[Proof of Theorem \ref{Th2}.]
	
        For $(\Leftarrow)$ of $(1)$, let $(a,b) \in \beta \circ^3 \gamma \cap \alpha \setminus \gamma$. Then there exist $x_1, x_2 \in A$ such that $a \ \beta\ x_1\ \gamma\ x_2\ \beta\ b$. Thus we have that the following relation holds:
		\begin{center}
			\begin{tikzcd}
				a \arrow[dd,"\b",dash]\arrow[r,"\a", dash]&b\arrow[dd,"\b",dash]
				\\
				\\
				x_1\arrow [r,"\g",dash]&x_2
			\end{tikzcd}
		\end{center}
		
		Hence, $((a,x_1), (b, x_2)) \in \alpha_0 \wedge \gamma_1$ and thus $\alpha_0 \wedge \gamma_1 > \gamma_0 \wedge \gamma_1$. Then, from Theorem \ref{Th1-2} and Lemma \ref{lemma1}, we have that the sublattice of $\con(\alg A(\beta))$ generated by $\{\alpha_0, \alpha_1, \gamma_0, \gamma_1, \beta_0\}$ is isomorphic to $\alg K$.

		For $(\Rightarrow)$ of $(1)$, let us suppose that the sublattice of $\con(\alg A(\beta))$ generated by $\{\alpha_0, \alpha_1, \gamma_0, \gamma_1, \beta_0\}$ is isomorphic to $\alg K$. Then, by Theorem \ref{Th1-2}, we have that $\alpha_0 \wedge \gamma_1$ and $\gamma_0 \wedge \alpha_1$ are not comparable and, by Lemma \ref{lemma1}, $\alpha_0 \wedge \gamma_1 > \gamma_0 \wedge \gamma_1$. Hence, there exists $((a_0, a_1), (b_0, b_1)) \in (\alpha_0 \wedge \gamma_1) \setminus (\gamma_0 \wedge \gamma_1)$. Thus we obtain that the following relations hold:
		
		\begin{center}
			\begin{tikzcd}
				a_0 \arrow[dd,"\b",dash]\arrow[r,"\a", dash]&b_0\arrow[dd,"\b",dash]
				\\
				\\
				a_1\arrow [r,"\g",dash]&b_1
			\end{tikzcd}
		\end{center}
		Thus $(a_0, b_0) \in (\beta \circ \gamma \circ \beta) \cap \alpha \setminus \gamma$ and this proves $(1)$.
		
		If  the sublattice of $\con(\alg A(\beta))$ generated by $\{\alpha_0, \alpha_1, \gamma_0, \gamma_1, \beta_0\}$ is not isomorphic to $\alg K$ then by Theorem \ref{Th1-2} is isomorphic to $\alg M_1$ and the thesis holds.
	\end{proof}
	
	The previous theorem emphasizes that the behaviour of $3$-permuting congruences has a deep impact on the structure of congruence lattices. Using this theorem we can provide examples of the two cases of Theorem \ref{Th1}. We can also observe that if $\gamma < \alpha$, then $(\beta \circ^3 \gamma) \cap \alpha \setminus \gamma \not= \emptyset$  if and only if $((\beta \circ^4 \gamma) \cap \alpha \setminus \gamma) \cup ((\gamma \circ^4 \beta)\cap \alpha \setminus \gamma) \not= \emptyset$. Thus, to find an example of an algebra $\alg A$ with $\alg M_1$ as the sublattice of $\con(\alg A(\beta))$ generated by $\{\alpha_0, \alpha_1, \gamma_0, \gamma_1, \beta_0\}$ we have to look at algebras whose $\beta$ and $\gamma$ at least do not $4$-permute.
	
	\begin{Ex}
		Let $\alg A$ be the algebra whose base set is $A = \{1,2,3,4,5,$ $6\}$,  with only projections as operations. Then the partitions $\beta = \mid 1, 2 \mid 3, 4 \mid 5, 6 \mid$, $\gamma = \mid 2, 3 \mid 4, 5 \mid$, and $\a = \mid 2, 3 \mid 4, 5 \mid 1, 6 \mid$  give $\mathbf{N}_5 \leq \con(\alg A)$ as in Figure \ref{N5}. Furthermore, using Theorem \ref{Th2}, we can see that $\alg M_1$ is the sublattice of $\con(\alg A(\beta))$ generated by $\{\alpha_0, \alpha_1, \gamma_0, \gamma_1, \beta_0\}$.
		
		On the other hand, let $\alg B$ be the algebra whose base set is $B = \{1,2,3,4\}$,  with only projections as operations. Then the partitions $\beta = \mid 1 ,3  \mid 2,4 \mid$ , $\gamma = \mid 1, 2\mid 3 \mid 4 \mid$, and $\a = \mid 1, 2\mid 3, 4 \mid$  give $\mathbf{N}_5 \leq \con(\alg B)$ as in Figure \ref{N5}. Furthermore, using Theorem \ref{Th2}, we can see that $\alg K$ is the sublattice of $\con(\alg B(\beta))$ generated by $\{\alpha_0, \alpha_1, \gamma_0, \gamma_1, \beta_0\}$.
	\end{Ex}

	\section{Duplications of Rods and Snakes}\label{Double}
	
	In this section we investigate the duplication of two particular families of lattices called \emph{rods} $\alg R_n$ and \emph{snakes} $\alg S_n$, $n \in \mathbb N$, $n \ge 1$ (see Figure \ref{rods}). These two families can be constructed using copies of $\alg M_3$ as building blocks and a formal definition of them could be the following. Let $(\alg M_3^n: n \in \mathbb N)$ be a countable family of copies of $\alg M_3$, where by $\a_i,\b_i,\g_i$ we denote the three atoms of the copy $\alg M_3^i$ from left to right, $\eta_i$ is the bottom and $\d_i$ is the top.  We glue the copies along a specified edge; this can be done since each upper edge of $\alg M_3$ is a lattice filter and each lower edge is a lattice ideal and they are isomorphic. Thus we define $\alg R_1 = \alg S_1 = \alg M_3$ and inductively
	\begin{enumerate}
		\ib to construct $\alg R_{n+1}$ we glue the filter $\{\g_n,\d_n\}$ of $\alg R_n$ to the ideal $\{\eta_{n+1},\a_{n+1}\}$ of $\alg M_3^{n+1}$ (Figure \ref{rods} right);
		\ib to construct $\alg S_{n+1}$, if $n$ is odd  we glue the filter $\{\g_n,\d_n\}$ of $\alg S_n$ to the ideal $\{\eta_{n+1},\a_{n+1}\}$ of $\alg M_3^{n+1}$, if $n$ is even  we glue the filter $\{\a_n,\d_n\}$ of $\alg S_n$ to the ideal $\{\eta_{n+1},\g_{n+1}\}$ of $\alg M_3^{n+1}$ (Figure \ref{rods} left).
	\end{enumerate}

	It is also clear that $\alg R_2 = \alg S_2 = \alg M_{3,3}$, see Figure \ref{ex:r2r2*}. These families have been introduced by A. Day in \cite{Day1975} together with $\alg R^*_n$ and $\alg S^*_n$, i.e. their counterparts with {\em the edges pulled apart}.
	
	\begin{figure}[htbp]
		\begin{center}
			\begin{tikzpicture}
				\draw (5,0) -- (5,2) -- (4,1) -- (5,0) -- (6,1) -- (5,2);
				\draw (6,1) -- (6,3) -- (5,2) -- (6,1) -- (7,2) -- (6,3);
				\draw (7.5,2.5) -- (7.5,4.5) -- (6.5,3.5) -- (7.5,2.5) -- (8.5,3.5) -- (7.5,4.5);
				\draw[dashed] (6,3) -- (6.5,3.5);
				\draw[dashed] (7,2) -- (7.5,2.5);
				\draw[fill] (5,0) circle [radius=0.05];
				\node[right] at (5,-0.1) {\footnotesize $\eta$};
				\draw[fill] (4,1) circle [radius=0.05];
				\draw[fill] (6,1) circle [radius=0.05];
				\node[left] at (4,1) {\footnotesize $\a_1$};
				\draw[fill] (5,1) circle [radius=0.05];
				\node[right] at (6,1) {\footnotesize $\g_1$};
				\draw[fill] (5,2) circle [radius=0.05];
				\node[left] at (5,2) {\footnotesize $\a_2$};
				\draw[fill] (5,1) circle [radius=0.05];
				\node[right] at (5,1) {\footnotesize $\b_1$};
				\draw[fill] (7,2) circle [radius=0.05];
				\node[right] at (7,2) {\footnotesize $\g_2$};
				\draw[fill] (6,3) circle [radius=0.05];
				\node[left] at (6,3) {\footnotesize $\a_3$};
				\draw[fill] (6,2) circle [radius=0.05];
				\node[right] at (6,2) {\footnotesize $\b_2$};
				\draw[fill] (7.5,2.5) circle [radius=0.05];
				\node[right] at (7.5,2.5) {\footnotesize $\g_{n-1}$};
				\draw[fill] (7.5,4.5) circle [radius=0.05];
				\node[right] at (7.5,4.6) {\footnotesize $\d$};
				\draw[fill] (6.5,3.5) circle [radius=0.05];
				\node[left] at (6.5,3.5) {\footnotesize $\a_n$};
				\draw[fill] (8.5,3.5) circle [radius=0.05];
				\node[right] at (8.5,3.5) {\footnotesize $\g_n$};
				\draw[fill] (7.5,3.5) circle [radius=0.05];
				\node[right] at (7.5,3.5) {\footnotesize $\b_n$};
				\node at (6,0.1) {\footnotesize $\alg R_n$};
				
				
				\draw (0,0) -- (2,2) -- (1,3) -- (1,1) -- (0,2) -- (0,0)-- (-1,1) -- (1,3);
				\draw (1,3) -- (0,4) -- (0,2) -- (-1,3) -- (0,4);
				\draw[dashed] (1,3) -- (2,4) -- (1,5) -- (0,4) -- (-1,5) -- (0,6);
				\draw (1,5) -- (2,6) -- (1,7) -- (1,5) -- (0,6) -- (1,7);
				\draw[fill] (0,0) circle [radius=0.05];
				\node[right] at (0,-0.1) {\footnotesize $\eta$};
				\draw[fill] (1,1) circle [radius=0.05];
				\node[right] at (1,1) {\footnotesize $\g_1$};
				\draw[fill] (1,3) circle [radius=0.05];
				\node[right] at (1,3) {\footnotesize $\g_3$};
				\draw[fill] (2,2) circle [radius=0.05];
				\node[right] at (2,2) {\footnotesize $\g_2$};
				\draw[fill] (0,2) circle [radius=0.05];
				\node[left] at (0,2) {\footnotesize $\a_2$};
				\draw[fill] (1,2) circle [radius=0.05];
				\node[right] at (1,2) {\footnotesize $\b_2$};
				\draw[fill] (-1,1) circle [radius=0.05];
				\node[left] at (-1,1) {\footnotesize $\a_1$};
				\draw[fill] (0,1) circle [radius=0.05];
				\node[right] at (0,1) {\footnotesize $\b_1$};
				\draw[fill] (0,3) circle [radius=0.05];
				\node[right] at (0,3) {\footnotesize $\b_3$};
				\draw[fill] (0,4) circle [radius=0.05];
				\draw[fill] (-1,3) circle [radius=0.05];
				\node[left] at (-1,3) {\footnotesize $\a_3$};
				\draw[fill] (1,5) circle [radius=0.05];
				\draw[fill] (2,6) circle [radius=0.05];
				\node[right] at (2,6) {\footnotesize $\g_n$};
				\draw[fill] (1,7) circle [radius=0.05];
				\node[right] at (1,7.1) {\footnotesize $\d$};
				\draw[fill] (0,6) circle [radius=0.05];
				\node[left] at (0,6) {\footnotesize $\a_n$};
				\draw[fill] (1,6) circle [radius=0.05];
				\node[right] at (1,6) {\footnotesize $\b_n$};
				\node at (1,0.1) {\footnotesize $\alg S_n$};
			\end{tikzpicture}
		\end{center}
		\caption{The families $\alg S_n$ and $\alg R_n$}\label{rods}
	\end{figure}

	In the same paper he showed that the $\alg R^*_n$ ($\alg S^*_n$), for $n\ge 1$, are all projective in the variety of modular lattices. Moreover for each $n$, $\alg R^*_n$ ($\alg S^*_n$) is a {\em projective cover}  of $\alg R_n$ ($\alg S_n$) (see \cite{Day1975} p. 154 for a definition). Since any lattice that is subdirectly irreducible and has a projective cover in a variety $\vv K$ of lattices is splitting in $\vv K$ (Theorem 3.7 in \cite{Day1975}) we conclude that each $\alg R_n$ ($\alg S_n$) is splitting in the variety of modular lattices.
	However $\alg R_n$ is splitting, but not projective and $\alg R^*_n$ is projective but not splitting (since it is not subdirectly irreducible).

        \begin{figure}[htbp]
		\begin{center}
			\begin{tikzpicture}
				\draw (4,0) -- (6.5,2.5) -- (5.5,3.5) -- (3,1) -- (4,0) -- (4,2) -- (5,1) -- (5.5,1.5) -- (5.5,3.5) -- (4.5,2.5) -- (5.5,1.5);
                \draw (-2,0) -- (0,2) -- (-1,3) -- (-3,1) -- (-2,0) -- (-2,2) -- (-1,3) -- (-1,1) -- (-2,2);
                 \draw[fill] (5.5,3.5) circle [radius=0.05];
                \node[right] at (5.5,3.6) {\footnotesize $\delta_2$};
                 \draw[fill] (4.5,2.5) circle [radius=0.05];
                \node[left] at (4.5,2.5) {\footnotesize $\alpha_2$};
                 \draw[fill] (4,2) circle [radius=0.05];
                \node[left] at (4,2) {\footnotesize $\delta_1$};
                 \draw[fill] (3,1.05) circle [radius=0.05];
                \node[left] at (3,1) {\footnotesize $\alpha_1$};
                
                \draw[fill] (5.5,2.5) circle [radius=0.05];
                \node[right] at (5.5,2.5) {\footnotesize $\beta_2$};
                \draw[fill] (4,1) circle [radius=0.05];
                \node[right] at (4,1) {\footnotesize $\beta_1$};
                 
                 \draw[fill] (6.5,2.5) circle [radius=0.05];
                \node[right] at (6.5,2.4) {\footnotesize $\gamma_2$};
                 \draw[fill] (5.5,1.5) circle [radius=0.05];
                \node[right] at (5.5,1.4) {\footnotesize $\eta_2$};
                 \draw[fill] (5,1) circle [radius=0.05];
                \node[right] at (5,0.9) {\footnotesize $\gamma_1$};
                 \draw[fill] (4,0) circle [radius=0.05];
                \node[right] at (4,-0.1) {\footnotesize $\eta_1$};
                \node[right] at (5,0.1) {\footnotesize $\alg M^*_{3,3}$};

                \draw[fill] (-1,2) circle [radius=0.05];
                \node[right] at (-1,2) {\footnotesize $\beta_2$};
                \draw[fill] (-2,1) circle [radius=0.05];
                \node[right] at (-2,1) {\footnotesize $\beta_1$};
                
                \draw[fill] (-1,3) circle [radius=0.05];
                \node[right] at (-1,3.1) {\footnotesize $\delta$};
                \draw[fill] (-2,2) circle [radius=0.05];
                \node[left] at (-2,2) {\footnotesize $\alpha_2$};
                \draw[fill] (-3,1) circle [radius=0.05];
                \node[left] at (-3,1) {\footnotesize $\alpha_1$}; 
                \draw[fill] (0,2) circle [radius=0.05];
                \node[right] at (0,1.9) {\footnotesize $\gamma_2$};
                \draw[fill] (-1,1) circle [radius=0.05];
                \node[right] at (-1,0.9) {\footnotesize $\gamma_1$};
                \draw[fill] (-2,0) circle [radius=0.05];
                \node[right] at (-2,-0.1) {\footnotesize $\eta$};
                \node[right] at (-1,0.1) {\footnotesize $\alg M_{3,3}$};
			\end{tikzpicture}
		\end{center}
		\caption{$\alg R_2 = \alg S_2 = \alg M_{3,3}$ and $\alg R^*_2 = \alg S^*_2 = \alg M^*_{3,3}$}\label{ex:r2r2*}
	\end{figure}
	
	We now state the aforementioned result of Freese and Lipparini;  we have to use the notion of {\em weak difference term} (see \cite{KearnesKiss2013}, Chapter $6$ for an extended treatment of the subject).
	
	\begin{Thm}[Freese-Lipparini, unpublished]\label{ralph}  Let $\vv V$ be a variety with a weak difference term and suppose that in $\vv V$ there is an algebra $\alg A$ with a nontrivial abelian congruence $\a$. Then there is field $\alg F$ such that for any $n$, the lattice of subspaces of an $n$-dimensional vector space over $\alg F$ is embeddable in $\op{Con}( \alg A^n(\a))$ (and the embedding is cover preserving).		
	\end{Thm}
	An application of the previous Theorem, combined with the fact that for any $n \in N$, $\alg R_n$ and $\alg S_n$ are characteristic free, yields the following corollary.
	
	\begin{Cor}\label{ralphcor} Let $\vv V$ be a variety having a weak difference term that is not congruence semidistributive; then $\alg R_n,\alg S_n \in \adm(\vv V)$ for all $n \in \mathbb N$.
	\end{Cor}
	The main result of the section (Theorem \ref{snakesThm}) generalizes Corollary \ref{ralphcor} since if a variety $\vv V$ has a weak difference term  and $\alg M_3 \leq \con(\alg A)$ for some $\alg A \in \vv V$, then the interval in $\con(\alg A)$  between bottom and top of the $\alg M_3$ consists of pairwise permuting congruences (see for instance \cite{KearnesKiss2013}). It follows that in such varieties, any $\alg M_3$ sublattice of a congruence lattice of an algebra in $\vv V$ satisfies the hypotheses of Theorem \ref{snakesThm}.
	
	In order to prove Theorem \ref{snakesThm} we will 
    apply the duplication technique to $\alg M_3$.
    From now on, until further notice, we will call the three atoms of it $\a,\b$ and $\g$ from the left to the right and without loss of generality we will apply the duplication construction to $\g$.
	
	\begin{Lem}\label{nozero} Let $\alg A$ be any algebra and suppose that $\alg M_3 \le \con(\alg A)$. Then the following are equivalent:
		\begin{enumerate}
			\item $\a_0 \meet \b_1$ and $\b_0 \meet \a_1$ are comparable;
			\item $\a_0 \meet \b_1 = \b_0 \meet \a_1 = 0_{\alg A(\g)}$.
		\end{enumerate}
	\end{Lem}
	\begin{proof} Suppose that $\a_0 \meet \b_1 \le \a_1 \meet \b_0$. Then
		$$
		\a_0 \meet \b_1 = \a_0 \meet \b_1 \meet \a_1\meet \b_0 = \eta_0 \meet \eta_1 = 0_{\alg A(\g)}.
		$$
		Moreover, let us consider the involutary automorphism which interchanges $\theta_0$ and $\theta_1$ for all $\theta \in \alg M_3$. Then if $\a_0 \meet \b_1 \le \a_1 \meet \b_0$ they must be equal and $(1) \Rightarrow (2)$. The converse is trivial.
	\end{proof}
	Let now $\alg A$ be any algebra such that $\alg M_3 \le \con(\alg A)$ and suppose that $\a,\b,\g \in \alg M_3$ $3$-permute (this is the same as asking that any pair of atoms $3$-permute). Let $(a,b) \in \a\setminus\b$. Then we have that
	$$
	(a,b) \in \a \join \b = \g \join \b = \g \circ \b\circ \g.
	$$
	Hence there are $c,d \in A$ with $a\mathrel{\g}c \mathrel{\b} d \mathrel{\g} b$. Then $((a,c), (b,d)) \in \a_0 \wedge \b_1$ with $a \not= b$ and thus $\alpha_0 \wedge \beta_1 \not= 0_{\alg A(\gamma)}$. By Lemma \ref{nozero}, we have that $\a_0 \meet \b_1$ and $\b_0 \meet \a_1$ are incomparable.
	
	Let us make a side remark. If $\vv V$ is a congruence $3$-permutable variety, then it is congruence modular and hence $\Con(\vv V)$ consists of modular lattices; so $\alg M_{3,3}$ (Figure \ref{fig14}) is splitting in $\Con(\vv V)$ and its splitting equation is (see \cite{Jonsson1968})
	$$
	x \meet (y \join (z \meet w)) \meet (z \join w) \le y \join (x\meet z) \join (x \meet w).
	$$
	If $\alg A \in \vv V$ and $\alg M_3 \le \con(\alg A)$, then the reader can check directly that the substitution
	$$
	x \mapsto \a_0 \meet \b_1 \ y \mapsto \a_0 \meet \a_1\ z \mapsto \b_0 \ w \mapsto \g_0
	$$
	falsifies the splitting equation, so $\alg M_{3,3} \in \Con(\vv V)$. This is nothing new, but it is interesting to observe that the argument works since we can evaluate the splitting equation by using the joins that we already know.
	
	However $3$-permutability is still not sufficient to prove what we want; we need that the congruences in $\alg M_3$ permute and in this case many congruences in
	$\op{Con}(\alg A(\g))$ are forced to permute as well, since the Correspondence Theorem \cite[Theorem $6.20$]{BurrisSanka}  preserves permutability for example. In the following lemma we list the permuting congruences needed to proceed in the proof of the main results of the section but the list is far from being exhaustive. We also evaluate some joins in $\op{Con}(\alg A(\g))$, using the following trivial observation: if $\lambda,\mu,\th$ are congruences such that
	$\lambda \join \mu \le \th$ and $\th \sse \lambda \circ \mu$, then $\lambda \join \mu = \th =  \lambda \circ \mu$. Thus $\lambda$ and $\mu$ permute.
	
	\begin{Lem}\label{permute} Let $\alg A$ be an algebra and suppose that $\alg M_3 \le \con(\alg A)$ consists of pairwise permuting congruences. Then the following equalities hold in $\op{Con}(\alg A(\g))$:
		\begin{enumerate}
			\item $\eta_0 \circ (\a_0 \meet \b_1) = \eta_0 \join (\a_0 \meet \b_1) = \a_0$;
            \item $\eta_0 \circ (\a_0 \meet \a_1) = \eta_0 \join (\a_0 \meet \a_1) = \a_0$;
			\item $(\a_0 \meet \b_1) \circ (\a_0 \meet \a_1) = (\a_0 \meet \b_1) \join (\a_0 \meet \a_1) = \a_0$;
			\item $\b_0 \circ \a_1 = \b_0 \join \a_1 = \d_0$.
		\end{enumerate}
	\end{Lem}
	\begin{proof}
		For $(1)$ let us prove $\a_0 \sse  \eta_0 \circ (\a_0 \meet \b_1)$. Thus let $((a,b),  (c,d)) \in\a_0$. Then $(a,c) \in \a$ and hence $(b,d) \in \g \circ \a\circ\g = \g \circ\a = \g \circ \b$.
		Thus there exists an $x \in A$ with $b\mathrel{\g} x \mathrel {\b} d$ and by transitivity $(a,x) \in \g$. Then:
		$$
		(a,b) \mathrel{\eta_0} (a,x)  \mathrel{\a_0 \meet\b_1} (c,d).
		$$
		Thus $\a_0 \sse  \eta_0 \circ (\a_0 \meet \b_1)$ and we can observe that $\eta_0 \circ (\a_0 \meet \b_1) \sse \eta_0 \join (\a_0 \meet \b_1) \sse \a_0$ and hence $(1)$ holds.
		\\
		The proof of $(2)$ is similar and we leave it to the reader.
		
		For $(3)$  let $(( a,b),  (c,d)) \in\a_0 = \eta_0 \circ (\a_0 \meet \a_1)$ by $(2)$. Then there exists $v \in A$ with:
		
		$$
		(a,b) \mathrel{\eta_0} (a,v)  \mathrel{\a_0 \meet\a_1}  (c,d).
		$$
		
		Thus $(a,v) \in \g$ and $(v,d) \in \a$. By transitivity $(b,v) \in \g \le \a \join \b = \b \circ \a$. Hence there exists an $x \in A$ with $b\mathrel{\b} x \mathrel{\a} v$; moreover, by transitivity, $(a,x) \in \g \circ \a = \a \circ \g$ so there exists $ y \in A$ such that $a \mathrel{\a} y\mathrel {\g} x$. Putting everything together we get
		$$
		(a,b) \mathrel{\a_0 \meet \b_1} (y,x) \mathrel{\a_0 \meet \a_1} (a,v) \mathrel{\a_0 \meet \a_1}  (c,d)
		$$
		and this proves that  $\alpha_0 = \eta_0 \circ (\a_0 \meet \a_1) \sse (\a_0 \meet \b_1) \circ (\a_0 \meet \a_1)$ and thus $(3)$. A graphical representation of this argument can be seen in Figure \ref{fig11}.
		
		\begin{figure}[htbp]
			\begin{center}
				\begin{tikzpicture}
					\node at (0,0) {$c$};
					\node at (1.2,0) {$d$};
					\node at (0,1.2) {$a$};
					\node at (1.2,1.2) {$v$};
					\node at (0,2.4) { $a$};
					\node at (1.2,2.4) { $b$};
					\draw (.2,0) -- (1,0);
					\draw (.2,1.2) -- (1,1.2);
					\draw (.2,2.4) -- (1,2.4);
					\draw (0,0.2) -- (0,1);
					\draw (0,1.4) -- (0,2.2);
					\draw (1.2,0.2) -- (1.2,1);
					\node at (2.5,1.2) {$\Longrightarrow$};
					\node at (3.8,0) {$c$};
					\node at (5,0) {$d$};
					\node at (3.8,1.2) {$a$};
					\node at (5,1.2) {$v$};
					\node at (3.8,2.4) {$a$};
					\node at (5,2.4) {$b$};
					\node at (4.5,1.8){$y$};
					\node at (5.7,1.8){$x$};
					\draw (4,0) -- (4.8,0);
					\draw (4,1.2) -- (4.8,1.2);
					\draw (4,2.4) -- (4.8,2.4);
					\draw (3.8,0.2) -- (3.8,1);
					\draw (3.8,1.4) -- (3.8,2.2);
					\draw (5,0.2) -- (5,1);
					\draw (3.9,1.3) -- (4.35,1.7);
					\draw (5.1,1.3) -- (5.55,1.7);
					\draw (4.7,1.8) -- (5.5,1.8);
					\draw (5.6,1.9) -- (5.15,2.3);
					\draw (3.9,2.3) -- (4.4,1.9);
					\node[above] at (0.6,0){\footnotesize $\g$};
					\node[above] at (0.6,1.2){\footnotesize $\g$};
					\node[above] at (0.6,2.4){\footnotesize $\g$};
					\node[left] at (0,0.6) {\footnotesize $\a$};
					\node[right] at (1.2,0.6) {\footnotesize $\a$};
					\node[above] at (4.4,0){\footnotesize $\g$};
					\node[below] at (4.4,1.2){\footnotesize $\g$};
					\node[above] at (4.4,2.4){\footnotesize $\g$};
					\node[left] at (3.8,0.6) {\footnotesize $\a$};
					\node[right] at (5,0.6) {\footnotesize $\a$};
					\node[above] at (5,1.8){\footnotesize $\g$};
					\node[left] at (4.55,1.5){\footnotesize $\a$};
					\node[right] at (5.3,1.5){\footnotesize $\a$};
					\node[right] at (5.35,2.2){\footnotesize $\b$};
					\node[right] at (4.20,2.2){\footnotesize $\a$};
				\end{tikzpicture}
			\end{center}
			\caption{}\label{fig11}
		\end{figure}
		For $(4)$ we can observe that $\eta_0 \circ \eta_1 = \gamma_0 = \gamma_1$. Thus:
		$$
		\beta_0 \circ \alpha_1 = (\beta_0 \circ \eta_0) \circ (\eta_1 \circ \alpha_1)  = \beta_0 \circ \gamma_0 \circ \gamma_1 \circ \alpha_1 = \delta_0 \circ \delta_1 = \delta_0
		$$
		and the thesis holds observing that $\b_0 \circ \a_1 \leq \b_0 \join \a_1 \leq \d_0$.
	\end{proof}
	Now we can show one of the main results of the section.
	
	\begin{proof} [Proof of Theorem \ref{M33}]
		First we can observe that $\op{Con}(\alg A(\g)/\eta_0)$ is isomorphic to $\con(\alg A)$ through the natural isomorphism $f(\lambda) = \lambda_0$. Thus the joins and meets in the part above $\eta_0$ in Figure \ref{fig14} are all respected. The first three claims of Lemma \ref{permute} ensure that all the joins of elements under $\alpha_0$ are respected. For the meets we can observe that $\eta_0 \wedge \alpha_0 \wedge \alpha_1 = \eta_0 \wedge \alpha_1 = \eta_0 \wedge \g_0 \meet \alpha_1 =\eta_0 \wedge \g_1 \meet \alpha_1 = \eta_0 \wedge \eta_1$. With the same argument we obtain that $\alpha_0 \wedge \b_1 \wedge \eta_0 = 0_{\alg A(\gamma)}$. A similar argument works also for $\alpha_0 \wedge \b_1 \wedge \alpha_0 \wedge \a_1$ $=$ $\alpha_0 \wedge \a_1 \wedge \b_0$ $ = $ $\alpha_0 \wedge \b_1 \wedge \b_0$ $=$ $\alpha_0 \wedge \a_1 \wedge \g_0$ $=$ $\alpha_0 \wedge \b_1 \wedge \g_0 = \eta_0 \wedge \eta_1$.
		
		For the remaining joins and the second claim we can see that, from Lemma \ref{permute}, the congruences under $\alpha_0$ pairwise permute. Moreover, since the congruences of the $\alg M_3$ pairwise permute, also the congruences above $\eta_0$ in Figure \ref{fig14} pairwise permute.
		
		Let us prove that the remaining pairs permute. We first make the trivial observation that if $\lambda \leq \nu$ then $\lambda \circ \nu = \nu$. We can use this, the relational products already computed, and Lemma \ref{permute} to prove that the following equalities hold:
		\begin{align*}
			&(\a_0 \wedge \a_1) \circ \beta_0 = (\a_0 \wedge \a_1) \circ \eta_0 \circ \beta_0 = \a_0 \circ \beta_0 = \delta_0
			\\&(\a_0 \wedge \b_1) \circ \beta_0 = (\a_0 \wedge \b_1) \circ \eta_0 \circ \beta_0 = \a_0 \circ \beta_0 = \delta_0
			\\&(\a_0 \wedge \a_1) \circ \g_0 = (\a_0 \wedge \a_1) \circ \eta_0 \circ \g_0 = \a_0 \circ \g_0 = \delta_0
			\\&(\a_0 \wedge \b_1) \circ \g_0 = (\a_0 \wedge \b_1) \circ \eta_0 \circ \g_0 = \a_0 \circ \g_0 = \delta_0.
		\end{align*}
		Thus all pairs of congruences in $$
		\{0_{\op{Con}(\alg A(\g))},\eta_0,\a_0 \meet\a_1,\a_0 \meet \b_1, \a_0,\b_0,\g_0,\d_0\}
		$$ permute.
	\end{proof}
	\begin{figure}[htbp]
		\begin{center}
			\begin{tikzpicture}
				\draw (0,0) -- (-1.5,1.5) -- (0,3) -- (0,0) -- (1.5,1.5) -- (0,3) -- (1.5,4.5) -- (3,3) -- (1.5,1.5) -- (1.5,4.5);
				\draw[fill] (0,0) circle [radius=0.05];
				\draw[fill] (-1.5,1.5) circle [radius=0.05];
				\draw[fill] (0,3) circle [radius=0.05];
				\draw[fill] (0,1.5) circle [radius=0.05];
				\draw[fill] (1.5,1.5) circle [radius=0.05];
				\draw[fill] (1.5,4.5) circle [radius=0.05];
				\draw[fill] (3,3) circle [radius=0.05];
				\draw[fill] (1.5,3) circle [radius=0.05];
				\node[right] at (0.1,-0.1) {\footnotesize $0_{\alg A(\g)}$};
				\node[left] at (-1.5,1.5) {\footnotesize $\a_0 \meet \a_1$};
				\node[left] at (0,3) {\footnotesize $\a_0$};
				\node[left] at (0,1.5) {\footnotesize $\a_0 \meet \b_1$};
				\node[right] at (1.5,1.5) {\footnotesize $\eta_0$};
				\node[right] at (1.6,4.6) {\footnotesize $\d_0$};
				\node[right] at (3,3) {\footnotesize $\g_0$};
				\node[right] at (1.5,3) {\footnotesize $\b_0$};
			\end{tikzpicture}
		\end{center}
		\caption{$\alg M_{3,3}$ as a sublattice of $\op{Con}(\alg A(\g))$}\label{fig14}
	\end{figure}
	Theorem \ref{M33} allows us to directly prove the main result of the section.
	
	\begin{proof} [Proof of Theorem \ref{snakesThm}]  Let  $\alg A$ be an algebra such that $\alg M_3 \le \con(\alg A)$ and all the congruences within $\alg M_3$ pairwise permute. We define by induction the sequence of algebras:
		\begin{align*}
			\alg A_1 &= \alg A
			\\ \alg A_i &= \alg A_{i-1}(\g_i)
		\end{align*}
		where $\g_2 = \g$ is an atom of $\alg M_3$ and $\g_i \in \con(\alg A_{i-1})$ is the kernel $\eta^{i-1}_0$ of the first projection from $\alg A_{i-1} \leq \alg A^2_{i-2}$ over the first factor, for all $i \geq 3$. We prove by induction on $n$ that $\alg R_n \leq \con(\alg A_n)$ and has permuting congruences, for all $n\ge 1$. Clearly, the claim holds for $n =1$ by hypothesis. 
		Let us consider $\con(\alg A_{n})$. By the inductive hypothesis we know that $\alg R_{n-1} \leq \con(\alg A_{n-1})$ and the congruences forming $\alg R_{n-1}$ pairwise permute. By the Correspondence Theorem,  the filter of  $\con(\alg A_{n})$ generated by $\eta^{n}_0$ is isomorphic to $\alg R_{n-1}$ and the congurences in this filter pairwise permute.  For this filter, we adopt the notation introduced in Figure \ref{rods}, with the addition of the superscript $n$ to specify that the congruences are in $\con(\alg A_{n})$ (see Figure \ref{Throds}). Furthermore, by Theorem \ref{M33}, the bottom $\alg M_{3}$ of  $\alg R_{n-1} \leq \con(\alg A_{n-1})$ produces a sublattice $\alg M_{3,3}$ of  $\con(\alg A_{n})$ with pairwise permuting congruences and where the atoms are $\{\alpha_n = (\a_{n-1})_0 \wedge  (\a_{n-1})_1, \b_n = (\a_{n-1})_0 \wedge  (\b_{n-1})_1, \g_n = \eta_0^n\}$ and $\a_1^n = (\a_{n-1})_0$ is the top element, see Figure \ref{Throds}. Thus $\alg R_{n} \leq \con(\alg A_{n})$ and $\alg R_{n} \in \adm(\vv V(\alg A))$. Hence, we have that all the congruences of the filter generated by $\eta_0^{n}$ and of the ideal generated by $\a^{n}_{1}$ pairwise permute.  Since $\a_n, \b_n \leq \a_i^{n}, \d^n$, for all $1 \leq i \leq n-1$ the remaining congruences we have to prove to permute are those in the set $\{\a_n, \b_n\}$ with the ones in the set $\{\b^{n}_i, \g^{n}_i \mid 2 \leq i \leq n-1\}$.  In particular, we have that $\b_n \circ \b^{n}_i = \b_n \circ \eta^{n}_0 \circ \b^{n}_i  = \a_1^{n} \circ \b^{n}_i = \a_{i+1}^{n}$, with $\a_{n}^{n} = \d^n$ and for all $1 \leq i \leq n-1$. The same proof also works for the other combinations of relational products involving elements from the set $\{\a_n, \b_n\}$ and elements of $\{\b^{n}_i, \g^{n}_i \mid 2 \leq i \leq n-1\}$. Hence, all the congruences of  $\alg R_{n} \leq \con(\alg A_{n})$ pairwise permute.
		Using the same strategy, we can prove by induction that $\alg S_n, \alg S^*_n, \alg R^*_n \leq \con(\alg A_n)$. In fact, the new $\alg M_3$ added at the bottom of the lattice can be glued either to the left or to the right and either tightly or loosely since $\alg R^*_{n-1} \leq \alg R_n$.
	\end{proof}
	
	\begin{figure}[htbp]
		\begin{center}
			\begin{tikzpicture}
				\draw (0,0) -- (-1.5,1.5) -- (0,3) -- (0,0) -- (1.5,1.5) -- (0,3) -- (1.5,4.5) -- (3,3) -- (1.5,1.5) -- (1.5,4.5);
				\draw (3.5,6.5)--(3.5,3.5)--(5,5)-- (3.5,6.5)-- (2,5) -- (3.5,3.5);
				\draw[dashed](2,5)--(1.5,4.5);
				\draw[dashed](3,3)--(3.5,3.5);
				\draw[fill] (0,0) circle [radius=0.05];
				\draw[fill] (-1.5,1.5) circle [radius=0.05];
				\draw[fill] (0,3) circle [radius=0.05];
				\draw[fill] (0,1.5) circle [radius=0.05];
				\draw[fill] (1.5,1.5) circle [radius=0.05];
				\draw[fill] (1.5,4.5) circle [radius=0.05];
				\draw[fill] (3,3) circle [radius=0.05];
				\draw[fill] (1.5,3) circle [radius=0.05];
				\draw[fill] (1.5,3) circle [radius=0.05];
				\draw[fill] (1.5,3) circle [radius=0.05];
				\draw[fill] (2,5) circle [radius=0.05];
				\draw[fill] (3.5,3.5) circle [radius=0.05];
				\draw[fill] (3.5,5) circle [radius=0.05];
				\draw[fill] (5,5) circle [radius=0.05];
				\draw[fill] (3.5,6.5) circle [radius=0.05];

				\node[right] at (0.1,-0.1) {\footnotesize $0_{\alg A(\g_{n-1})}$};
				\node[left] at (-1.5,1.5) {\footnotesize $\alpha_n$};
				\node[left] at (0,3) {\footnotesize $\a_1^{n}$};
				\node[left] at (0,1.5) {\footnotesize $\b_n$};
				\node[right] at (1.5,1.5) {\footnotesize $\eta^{n}_0 = \g_n$};
				\node[right] at (0.8,4.6) {\footnotesize $\a_2^{n}$};
				\node[right] at (1,5.1) {\footnotesize $\a_{n-1}^{n}$};
				\node[right] at (3,3) {\footnotesize $\g^{n}_1$};
				\node[right] at (3.5,3.5) {\footnotesize $\g^{n}_{n-2}$};
				\node[right] at (5,5) {\footnotesize $\g^{n}_{n-1}$};
				\node[right] at (3.5,5) {\footnotesize $\b^{n}_{n-1}$};
				\node[right] at (3.5,6.6) {\footnotesize $\delta^{n}$};
				\node[right] at (1.5,3) {\footnotesize $\b^{n}_1$};
			\end{tikzpicture}
		\end{center}
		\caption{$\alg R_{n}$ as a sublattice of $\con(\alg A_{n})$}\label{Throds}
	\end{figure}
	It is worth noting that the proof of Theorem \ref{snakesThm} also works for \emph{generalized snakes}, i.e. lattices $\alg G_n$ obtained by glueing $n$ copies of $\alg M_3$ (with atoms $\alpha,\b,\gamma$) either by collapsing the top-right edge $(\gamma, 1)$ of one instance with the bottom-left $(0, \a)$ of the other copy (as seen in the rods) or by merging the top-left edge $(\a, 1)$ with the  bottom-right $(0, \g)$ (as seen in the even steps of snakes).

	\section{Conclusions and problems}\label{Conclusions}
	
	Let $\vv V$ be a variety; if $\vv V$ has a weak difference term  and $\alg M_3 \leq \con(\alg A)$ for some $\alg A \in \vv V$, then the interval in $\con(\alg A)$  between bottom and top of the $\alg M_3$ consists of pairwise permuting congruences (see for instance \cite{KearnesKiss2013}). It follows that in such varieties, any $\alg M_3$ as a sublattice of a congruence lattice of an algebra in $\vv V$ satisfies the hypotheses of Theorem \ref{snakesThm}.
	
	However, in order to apply the results in Section \ref{Double} we do not need that {\em all} the  $\alg M_3$'s in all algebras in a variety consist of permuting congruences (as in case of a week difference term): we only need one of those instances in $\vv V$ to show that $\alg R_n, \alg S_n \in \adm(\vv V)$ for all $n$. Now K. Kearnes (in a private communication) showed us that given a variety $\vv V$ that is not congruence meet semidistributive there is  an $\alg A \in \vv V$ such that $\alg M_3 \le \con(\alg A)$ and the congruences in the $\alg M_3$ permute.  This is connected to a class of problems coming from the concept of {\em omission class}. Let $\Gamma$ be a class of lattices; then the \emph{omission class of $\Gamma$} is
$$
	\mathfrak F(\Gamma)=\{\vv V: \alg L \notin \adm(\vv V), \text{for all}\ \alg L \in \Gamma\}.
$$
Kearnes' result above implies that
 $$
 \mathfrak F(\alg M_3) = \bigcup_{n \ge 1}\mathfrak F(\alg R_n) = \bigcup_{n \ge 1}\mathfrak F(\alg S_n).
 $$
 This  is of interest since the omission class of $\alg M_3$ coincides with the class of meet semidistributive varieties (see below).
	
 There are also other natural questions. As far as we know there are only   four nontrivial Mal'cev classes  that are equal to the omission class of $\alg L$ for some finite subdirectly irreducible lattice $\alg L$:
	\begin{enumerate}
		\ib the class $T$ of all varieties satisfying a nontrivial idempotent Mal'cev condition ($\mathfrak{F}(\alg D_1)$, where $\alg D_1$ is displayed in Figure \ref{M1} as the filter generated by $\g_0 \wedge \g_1$ \cite[Theorem 4.23]{KearnesKiss2013});
		\ib the class $L$ of all varieties $\vv V$ such that $\Con(\vv V)$ is a proper subvariety of the variety  of lattices ($\mathfrak{F}(\alg D_2)$, where $\alg D_2$ is the dual of $\alg D_1$ \cite[Theorem 8.11]{KearnesKiss2013}) ;
		\ib the class $M$ of all congruence modular varieties ($\mathfrak{F}(\alg N_5)$ from \cite{Dedekind1900});
		\ib the class $SD_\meet$ of all congruence meet semidistributive varieties ($\mathfrak{F}(\alg M_3)$ by \cite[Theorem 8.1]{KearnesKiss2013}).
	\end{enumerate}
	So it makes sense to ask if these are the only Mal'cev classes with this property. A harder question is: can we give necessary and/or sufficient condition on a finite subdirectly irreducible lattice $\alg L$, so that the omission class of $\alg L$ is a Mal'cev class? We will deal with this set of problems in
a further paper.

	\section*{Acknowledgements}
	The authors thank the referee and editor for their useful suggestions and for the discussion raised about the conclusions of the paper.

\providecommand{\bysame}{\leavevmode\hbox to3em{\hrulefill}\thinspace}
\providecommand{\MR}{\relax\ifhmode\unskip\space\fi MR }
\providecommand{\MRhref}[2]{%
  \href{http://www.ams.org/mathscinet-getitem?mr=#1}{#2}
}
\providecommand{\href}[2]{#2}

\end{document}